\documentclass[12pt]{article}

\newlength{\stefan}
\usepackage{amsmath, amsthm, amsfonts, makeidx}
\usepackage{multind}
\DeclareMathSymbol{\subsetneq}{\mathord}{AMSb}{"26}

\newtheorem{lemma}{Lemma}[section]

\newtheorem{proposition}[lemma]{Proposition}

\theoremstyle{definition}

\newcommand{\mb}{\mathbb}

\newcommand{\C}{\mb{C}}

\newcommand{\N}{\mb{N}}

\title{Infinitely generated Derksen and Makar-Limanov invariant\footnote{AMS classification 13N15,13A50,14R20}}
\author{Stefan Maubach\footnote{Funded by a Veni-grant of the council for the
physical sciences, Netherlands Organisation for scientific research (NWO)}}

\begin{document}

\maketitle

\begin{abstract}
In this paper, we give an example of a finitely generated 3-dimensional $\C$-algebra which  has infinitely generated Derksen invariant as well as Makar-Limaonv invariant.
\end{abstract}

\section{Introduction and tools}

The Derksen invariant and Makar-Limanov invariant are useful tools to recognize if two varieties or rings are not isomorphic. 
Both invariants use locally nilpotent derivations: if $A$ is a commutative $k$-algebra (where $k$ is a field of characteristic zero), 
then $D$ is a derivation if $D$ is $k$-linear and satisfies the Leibniz rule: $D(ab)=aD(b)+bD(a)$. A derivation is locally nilpotent if 
for each $a\in A$ we can find some $n\in \N$ such that $D^n(a)=0$. The kernel of a derivation, denoted by $A^D$, is the set of all elements that are mapped 
to zero under the derivation $D$. The Makar-Limanov invariant is defined
as the intersection of all kernels of locally nilpotent derivations, while the Derksen invariant is defined 
as the smallest algebra containing the kernels of all nonzero locally nilpotent derivations.

In the paper \cite{R++} the question was posed if the Derksen invariant could be infinitely generated.
In this paper we give an example of an infinitely generated Derksen invariant of a finitely generated $\C$-algebra.
It will be at the same time an example of an infinitely generated Makar-Limanov invariant, as in this example, the Derksen invariant is equal to the Makar-Limanov invariant.
By now, there are many examples of cases of ``nice'' subrings that are not finitely generated \cite{BD97,DF,Mau,Kur1,Kur2}. In regard of this, the author would like to remark
that it will pay off to consider theorems as general as possible (with respect to not restricting to finitely generated algebras).

{\sl Notations:} If $R$ is a ring, then $R^{[n]}$ denotes the polynomial ring in $n$ variables over $R$.
We will use the letter $k$ for a field of characteristic zero, and $K$ for its algebraic closure.
Denote by $\partial_x$ the derivative with respect to $x$. By $LND(A)$ we will denote the set of all locally nilpotent derivations on a ring $A$. 

The technique used for constructing this example, is based on the following general  idea: grab a locally nilpotent derivation $D$ on a polynomial ring $K^{[n]}$ having the required
properties (non-finitely generated kernel).
Then construct an equation $f$ (or several equations $f_i$) which forces $A:=K^{[n]}/(f)$ (or $A:=K^{[n]}/(f_1,\ldots,f_m)$) 
to be a ring that has only one locally nilpotent derivation (up to multiplication with an element of $A$), namely $D \mod (f_1,\ldots,f_m)$. 
Hopefully it still has the required properties (infinitely generated kernel).

Well-known facts that we need  are the following:

\begin{lemma}\label{extra}
Let $D\in LND(A)$ where $A$ is a domain. \\
(1) Then $D(A^*)=0$.\\
(2) If $D(ab)=0$ with $a,b$ both nonzero, then $D(a)=D(b)=0$.\\
(3) If $fD\in LND(A)$ where $f\not = 0$ then $D(f)=0$ and $D\in LND(A)$.
\end{lemma}

\section{The  example}

This example is inspired by the example of Bhatwadekar and Dutta in \cite{BD97}.
We will write small letters for capital letters modulo a relation: for example, $a$ below is defined as $A+(A^3-B^2)$. 

Define $R:=\C[A,B]/(A^3-B^2)=\C[a,b]\cong\C[T^2,T^3]\subseteq \C[T]$.
Define $S:=R[X,Y,Z]/(Z^2-a^2(aX+bY)^2-1)=R[X,Y][z]$.
We leave it to the reader to check that $S$ is a domain.
We will first try to find all locally nilpotent derivations on this ring.

\begin{lemma}\label{E1.L1}
Let $D\in LND(S)$. Then $D(a)=D(b)=D(aX+bY)=D(z)=0$.
\end{lemma}

\begin{proof}
Since $(z-a^2X-abY)(z+a^2X+abY)=1$, we have by lemma \ref{extra} part 1 that $D(z-a^2X-abY)=D(z+a^2X+abY)=0$, and thus since $D$ is $\C$-linear,
$D(z)=D(a^2X+abY)=0$.
Because of lemma \ref{extra} part 2, we have $D(a)=D(aX+bY)=0$. Since $0=3a^2D(a)=D(a^3)=D(b^2)=2bD(b)$ we have $D(b)=0$.
\end{proof}

\begin{lemma}\label{E1.L2}
$LND(S)=S^D \cdot D$ where $D:= b\partial_X-a\partial_Y$.
\end{lemma}

\begin{proof}
Let $D\in LND(S)$. Then $aD(X)=-bD(Y)$ by lemma \ref{E1.L1}. Seeing $S$ as a subring of $B:=C\oplus C\bar{z}=C[Z]/(Z^2-T^8(X+TY)^2-1)$ where $C:=\C[T][X,Y]$, we write
$D(X)=f_0+\bar{z}f_1, D(Y)=g_0+\bar{z}g_1$ where $f_0,f_1, g_0,g_1\in C$. Since $D(X)=-TD(Y)$ we see that $T$ divides $f_0$ and $f_1$. But since both are in
$\C[T^2,T^3,X,Y]$ we know that even $T^2$ divides them. So, $D(X)=T^2h_0+T^2\bar{z}h_1$ where $h_i:=f_iT^{-2}$. But now $T(h_0+h_1\bar{z})=-g_0-g_1\bar{z}$, so $T$ divides the $g_i$.
Again, we have that even $T^2$ divides $g_i$, which then gives that $T$ divides $h_i$. In the end, $T^3$ divides $D(X)$ and $T^2$ divides $D(Y)$.
Write $D(X)=T^3f=T^3f_0+T^3zf_1, D(Y)=T^2g=T^2g_0+T^2zg_1$. Then $T^3f=D(X)=-TD(Y)=-T^3g$ so $f_0=-g_0,f_1=-g_1,$ and thus $f=-g$, and $D(X)=bf, D(Y)=af$. 
Since $D(z)=D(a)=D(b)=0$, $D=f( b\partial_X-a\partial_Y)$.
By lemma \ref{extra} we have that $f\in S^D$ and so we are done.
\end{proof}

\begin{proposition}\label{E1.P3}
Let $D:= T^3\partial_X-T^2\partial_Y$ on $S$ as before. Then  $S^D$ is not finitely generated as a $\C$-algebra.
\end{proposition}

\begin{proof}
Examining the natural extension of $D$ on $B:=\C[T][X,Y][z]$ it is easy to determine that $B^D=\C[T,z,X+TY]$.
Now $S^D=B^D\cap S$, as can be easily checked. Defining $P:=X+TY$, we are done if we show that $\C[T,z,P]\cap\C[T^2,T^3,z,X,Y]$ is not finitely generated.
We will do this by writing elements in a unique way in a representant system. \\
{\sc claim:} If $F\in S^D\backslash \C[T,z]$, then $F\in (T^2,T^3)S$. \\
{\sc proof of claim:} $S^D\mod(T^2,T^3)=\C[T,X+TY,z]/(T^2,z^2-1) \cap \C[z,X,Y]/(z^2-1)$. 
If $F\in S^D$ is nonzero, then 
\[ F\mod(T^2)=\sum_{i=0}^n f_i(\bar{T},\bar{z})(X+\bar{T}Y)^i =\sum_{i=0}^n f_i(\bar{T},\bar{z})(X^i+i\bar{T}X^{i-1}Y) \]
where $f_n\not = 0$. But since $F\mod(T^2,T^3)S \in \C[\bar{z},X,Y]$, this implies $n=0$. 
So $F\in \C[T,z]+(T^2,T^3)S$, which proves the claim.

If $S^D$ is finitely generated, then $S^D=\C[T^2,T^3,z,F_1,\ldots,F_n]$ where $F_i\in (T^2)B$ (by the claim).
Let $d$ be the maximum of the $X,Y$-degree of the $F_i$.
Take $T^2P^{d+1}$, which is in $S^D$.
Write $T^2P^{d+1}=f(F_1,\ldots,F_n)$ where $f$ has coefficients in $\C[T^2,T^3]\oplus
\C[T^2,T^3]z$. Computing modulo $(T^4)B$ (or $(T^4,T^5,T^6)S$) we see that
$T^2P^{d+1}=f_0+f_1F_1+\ldots+f_nF_n \mod{T^4}$ where $f_i\in \C[\bar{T}^2,\bar{T}^3]\oplus
\C[\bar{T}^2,\bar{T}^3]\bar{z}$, as $F_iF_j\in (T^4)B$ for all $i,j$. 
Now comparing the $X,Y$-degree from the right hand side (which is at most $d$) to the 
$X,Y$ degree  on the left hand side (which is $d+1$) we get a contradiction, showing that the assumption ``$S^D$ is finitely generated'' is wrong.
\end{proof}

{\bf Remark:}
The example in this paper is no UFD.
In the paper \cite{FM} an example of a $\C$-algebra UFD of dimension 6 is given, which has infinitely generated Derksen and Makar-Limanov invariant.

\noindent
{\sc Author's address}\\
\small
Stefan Maubach\\
Radboud University Nijmegen\\
\small Toernooiveld 1, The Netherlands\\ 
\small s.maubach@science.ru.nl

\end{document}